\documentclass[12pt,reqno]{article}
\usepackage{mathtools}
\usepackage[usestackEOL]{stackengine}[2013-10-15]
\usepackage[usenames]{color}
\usepackage{amsmath,amsthm,amssymb,tikz,enumerate}
\usetikzlibrary{decorations.pathreplacing}
\usepackage{color}
\usepackage{fullpage}
\usepackage{float}
\usepackage{epsf}
\newcommand\tab[1][1cm]{\hspace*{#1}}
\usepackage[colorlinks=true,linkcolor=webgreen,filecolor=webbrown,citecolor=webgreen]{hyperref}
\usepackage{graphicx}

\definecolor{webgreen}{rgb}{0,.5,0}
\definecolor{webbrown}{rgb}{.6,0,0}

\begin{document}

\begin{center}
\end{center}
\theoremstyle{plain}
\newtheorem{theorem}{Theorem}
\newtheorem{conclusion}{Conclusion}
\newtheorem{corollary}{Corollary}
\newtheorem{lemma}{Lemma}
\theoremstyle{definition}
\newtheorem{definition}{Definition}
\begin{center}
\vskip 1cm{\LARGE\bf The Relationship Between Pascal's Triangle and Random Walks}
\vskip 1cm
\large
Tonia Bell, Shakuan Frankson, Nikita Sachdeva, and Myka Terry\\ 
Mathematics Department\\
SPIRAL Program at Morgan State University, Baltimore, MD\\
\href{mailto:tb8755a@student.american.edu}{\tt tb8755a@student.american.edu}\\
\href{mailto:shakuan.frankson@bison.howard.edu}{\tt shakuan.frankson@bison.howard.edu}\\
\href{mailto:nikita.sachdeva@rutgers.edu}{\tt nikita.sachdeva@rutgers.edu}\\
\href{myter1@morgan.edu}{\tt myter1@morgan.edu}\\

\end{center}
\vskip .2 in

\begin{abstract}
\noindent Random walks are a series of up, down, and level steps that enumerate distinct paths from $(0,0)$ to $(2n,0)$, where $n$ is the semi-length of the path. We used these paths to analyze Catalan, Schr\"{o}der, and Motzkin number sequences through a combination of matrix operations, quadratic functions, and inductive reasoning. Our results revealed a number of distinct patterns, some unnamed, between these number sequences and Pascal's triangle that can be explained through generating functions, first returns, group theory, and the Riordan matrix.
\\
\\
Various proofs and properties of these number sequences are provided, including each generating function, their respective first returns, and matrix properties. These findings lead to a deeper understanding of combinatorics and graph theory.
\end{abstract}
\emph{\section{Introduction}}
\noindent In the Cartesian plane, a single vector or combination of vectors can create a path, where each vector represents a \emph{step} towards a unique xy-coordinate. In this paper, we restricted the movement of the vectors to the first quadrant. Using the $(1,1)$ vector as the up step (due to its positive slope) and $(1,-1)$ vector as the down step (due to its negative slope), along with vectors $(1,0)$ and $(2,0)$ as the horizontal steps (due to their slopes of 0), we generated the Catalan, Schr\"{o}der, and Motzkin numbers along the x-axis.\\

\noindent The  original  Dyck  path uses steps along the vectors  $(1,1)$  and  $(1,-1)$  and remains above the $x$-axis, extending from $(0,0)$ to $(2n,0)$, where $n$ is the semi-length. The path is also considered to be the Catalan path since it enumerates the Catalan numbers on the bottom row along the x-axis. Schr\"{o}der  numbers follow the same pattern as the Catalan numbers, while also having the level step, $(2,0)$. Their path is referred to as a Schr\"{o}der path. The Motzkin numbers of length $n$ have the smaller level step of $(1,0)$ and remain above the $x$ -axis from $(0,0)$ to $(n,0)$. Similarly, these paths are referred to as Motzkin paths.\\

\noindent Pascal's Triangle will be referenced to in Section 3 since the method for finding the binomial coefficients will parallel the process of constructing the Catalan, Schr\"{o}der, and Motzkin paths. In this case, they also share a generalized recursive formula that determines the next row of values, or sequence, when two or more rows are multiplied together.
\\

\subsection{Catalan Numbers}

Catalan numbers, denoted $C(x)$ or $C$, are a variation of Dyck paths that start at $(0,0)$ and end at $(2n,0)$. Recall the following vectors as the only two found in Catalan paths:
\begin{center}
\begin{tikzpicture}
	\draw[->](3,0) -- (3.8,-0.7);
	\draw[->](7.2,-0.7) -- (8,0);
	
\end{tikzpicture}		

\tab[0.00015cm]$(1,-1)$\tab[3cm]$(1,1)$\\ 
\end{center}
which can be seen in the following lattice: 
\begin{center}

\begin{tikzpicture}[ultra thick]
\draw [help lines] (0,0) grid (11,6);
\coordinate (aux) at (0,0);
\draw[->] [color=gray] (0,0)--(0,6);
\draw[->] [color=gray] (0,0)--(11,0);
\draw[->] [color=teal] (0,0)--(5,5);
\draw[->] [color=teal] (5,5)--(10,0);
\draw[->] [color=teal] (1,1)--(2,0);
\draw[->] [color=teal] (2,0)--(3,1);
\draw[->] [color=teal] (3,1)--(4,0);
\draw[->] [color=teal] (4,0)--(5,1);
\draw[->] [color=teal] (5,1)--(6,0);
\draw[->] [color=teal] (6,0)--(7,1);
\draw[->] [color=teal] (7,1)--(8,0);
\draw[->] [color=teal] (8,0)--(9,1);
\draw[->] [color=teal] (9,1)--(10,0);
\draw[->] [color=teal] (2,2)--(3,1);
\draw[->] [color=teal] (3,1)--(4,2);
\draw[->] [color=teal] (4,2)--(5,1);
\draw[->] [color=teal] (5,1)--(6,2);
\draw[->] [color=teal] (6,2)--(7,1);
\draw[->] [color=teal] (7,1)--(8,2);
\draw[->] [color=teal] (3,3)--(4,2);
\draw[->] [color=teal] (4,2)--(5,3);
\draw[->] [color=teal] (5,3)--(6,2);
\draw[->] [color=teal] (6,2)--(7,3);
\draw[->] [color=teal] (4,4)--(5,3);
\draw[->] [color=teal] (5,3)--(6,4);

\node at (1,1) [circle,draw=purple!100,fill=purple!0,minimum size=0.95cm] {};
\node at (3,1) [circle,draw=purple!100,fill=purple!0,minimum size=0.95cm] {};
\node at (5,1) [circle,draw=purple!100,fill=purple!0,minimum size=0.95cm] {};
\node at (7,1) [circle,draw=purple!100,fill=purple!0,minimum size=0.95cm] {};
\node at (9,1) [circle,draw=purple!100,fill=purple!0,minimum size=0.95cm] {};
\node at (2,2) [circle,draw=purple!100,fill=purple!0,minimum size=0.95cm] {};
\node at (4,2) [circle,draw=purple!100,fill=purple!0,minimum size=0.95cm] {};
\node at (6,2) [circle,draw=purple!100,fill=purple!0,minimum size=0.95cm] {};
\node at (8,2) [circle,draw=purple!100,fill=purple!0,minimum size=0.95cm] {};
\node at (3,3) [circle,draw=purple!100,fill=purple!0,minimum size=0.95cm] {};
\node at (5,3) [circle,draw=purple!100,fill=purple!0,minimum size=0.95cm] {};
\node at (7,3) [circle,draw=purple!100,fill=purple!0,minimum size=0.95cm] {};
\node at (4,4) [circle,draw=purple!100,fill=purple!0,minimum size=0.95cm] {};
\node at (6,4) [circle,draw=purple!100,fill=purple!0,minimum size=0.95cm] {};
\node at (5,5) [circle,draw=purple!100,fill=purple!0,minimum size=0.95cm] {};
\node at (0,0) [circle,draw=purple!100,fill=purple!0,minimum size=0.95cm] {};
\node at (2,0) [circle,draw=purple!100,fill=purple!0,minimum size=0.95cm] {};
\node at (4,0) [circle,draw=purple!100,fill=purple!0,minimum size=0.95cm] {};
\node at (6,0) [circle,draw=purple!100,fill=purple!0,minimum size=0.95cm] {};
\node at (8,0) [circle,draw=purple!100,fill=purple!0,minimum size=0.95cm] {};
\node at (10,0) [circle,draw=purple!100,fill=purple!0,minimum size=0.95cm] {};

\node (a) at (1,1) {1};
\node (b) at (3,1) {2};
\node (c) at (5,1) {5};
\node (d) at (7,1) {14};
\node (e) at (9,1) {42};
\node (f) at (2,2) {1};
\node (g) at (4,2) {3};
\node (h) at (6,2) {9};
\node (i) at (8,2) {28};
\node (j) at (3,3) {1};
\node (k) at (5,3) {4};
\node (l) at (7,3) {14};
\node (m) at (4,4) {1};
\node (n) at (6,4) {5};
\node (o) at (5,5) {1};
\node (p) at (0,0) {1};
\node (q) at (2,0) {1};
\node (r) at (4,0) {2};
\node (s) at (6,0) {5};
\node (t) at (8,0) {14};
\node (u) at (10,0) {42};

\end{tikzpicture}
\begin{figure}[h!]
\caption{Catalan Path, C(x)}
\end{figure}
\end{center}

\noindent The following is the generating function for the Catalan numbers and will be proved in Section 2.
\begin{equation}
C(x)=\frac{1-\sqrt{1-4x}}{2x} \label{eq:1}
\end{equation}
\noindent This produces the infinite Catalan sequence, which can be found along the $x$-axis of the lattice in Figure 2:\\
\begin{equation}
C(x)=1+x+2x^2+5x^3+14x^4+42x^5+...\\
\end{equation}
\subsection{Schr\"{o}der Numbers}
There are two types of Schr\"{o}der numbers seen in combinatorics: large Schr\"{o}der numbers denoted, $S(x)$ or $S$, and small Schr\"{o}der numbers denoted, $s(x)$ or $s$. They are composed of the vectors:
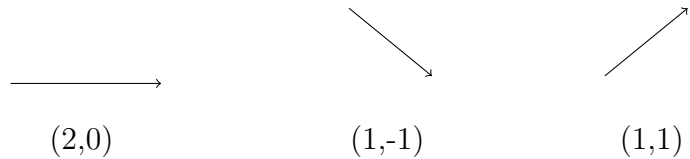
\begin{center}
\begin{tikzpicture}
	\draw[->](0,-1) -- (2,-1);
	\draw[->](4.5,0) -- (5.6,-0.9);
	\draw[->](7.9,-0.9) -- (9,0);
	
\end{tikzpicture}		
\end{center}

\tab[3.65cm](2,0)\tab[3.15cm](1,-1)\tab[2.6cm](1,1)\\

\noindent These vectors create unique Schr\"{o}der  paths for the large Schr\"{o}der and small Schr\"{o}der numbers, respectively. The sole distinction between them is $s(x)$ does not consist of the level step, $(2,0)$, along the x-axis, while $S(x)$ does. Below is the lattice for both series. 
\\

\begin{center}
\begin{tikzpicture}[ultra thick]
\draw [help lines] (0,0) grid (11,6);
\draw[->] [color=gray] (0,0)--(0,6);
\draw[->] [color=gray] (0,0)--(11,0);
\coordinate (aux) at (0,0);
\draw[->] [color=teal] (0,0)--(5,5);
\draw[->] [color=teal] (5,5)--(10,0);
\draw[->] [color=teal] (1,1)--(2,0);
\draw[->] [color=teal] (2,0)--(3,1);
\draw[->] [color=teal] (3,1)--(4,0);
\draw[->] [color=teal] (4,0)--(5,1);
\draw[->] [color=teal] (5,1)--(6,0);
\draw[->] [color=teal] (6,0)--(7,1);
\draw[->] [color=teal] (7,1)--(8,0);
\draw[->] [color=teal] (8,0)--(9,1);
\draw[->] [color=teal] (9,1)--(10,0);
\draw[->] [color=teal] (2,2)--(3,1);
\draw[->] [color=teal] (3,1)--(4,2);
\draw[->] [color=teal] (4,2)--(5,1);
\draw[->] [color=teal] (5,1)--(6,2);
\draw[->] [color=teal] (6,2)--(7,1);
\draw[->] [color=teal] (7,1)--(8,2);
\draw[->] [color=teal] (3,3)--(4,2);
\draw[->] [color=teal] (4,2)--(5,3);
\draw[->] [color=teal] (5,3)--(6,2);
\draw[->] [color=teal] (6,2)--(7,3);
\draw[->] [color=teal] (4,4)--(5,3);
\draw[->] [color=teal] (5,3)--(6,4);
\draw[->] [color=teal] (2,0)--(4,0);
\draw[->] [color=teal] (4,0)--(6,0);
\draw[->] [color=teal] (6,0)--(8,0);
\draw[->] [color=teal] (8,0)--(10,0);
\draw[->] [color=teal] (0,0)--(10,0);
\draw[->] [color=teal] (1,1)--(9,1);
\draw[->] [color=teal] (2,2)--(8,2);
\draw[->] [color=teal] (3,3)--(7,3);
\draw[->] [color=teal] (4,4)--(6,4);

\node at (1,1) [circle,draw=purple!100,fill=purple!0,minimum size=0.95cm] {};
\node at (3,1) [circle,draw=purple!100,fill=purple!0,minimum size=0.95cm] {};
\node at (5,1) [circle,draw=purple!100,fill=purple!0,minimum size=0.95cm] {};
\node at (7,1) [circle,draw=purple!100,fill=purple!0,minimum size=0.95cm] {};
\node at (9,1) [circle,draw=purple!100,fill=purple!0,minimum size=0.95cm] {};
\node at (2,2) [circle,draw=purple!100,fill=purple!0,minimum size=0.95cm] {};
\node at (4,2) [circle,draw=purple!100,fill=purple!0,minimum size=0.95cm] {};
\node at (6,2) [circle,draw=purple!100,fill=purple!0,minimum size=0.95cm] {};
\node at (8,2) [circle,draw=purple!100,fill=purple!0,minimum size=0.95cm] {};
\node at (3,3) [circle,draw=purple!100,fill=purple!0,minimum size=0.95cm] {};
\node at (5,3) [circle,draw=purple!100,fill=purple!0,minimum size=0.95cm] {};
\node at (7,3) [circle,draw=purple!100,fill=purple!0,minimum size=0.95cm] {};
\node at (4,4) [circle,draw=purple!100,fill=purple!0,minimum size=0.95cm] {};
\node at (6,4) [circle,draw=purple!100,fill=purple!0,minimum size=0.95cm] {};
\node at (5,5) [circle,draw=purple!100,fill=purple!0,minimum size=0.95cm] {};
\node at (0,0) [circle,draw=purple!100,fill=purple!0,minimum size=0.95cm] {};
\node at (2,0) [circle,draw=purple!100,fill=purple!0,minimum size=0.95cm] {};
\node at (4,0) [circle,draw=purple!100,fill=purple!0,minimum size=0.95cm] {};
\node at (6,0) [circle,draw=purple!100,fill=purple!0,minimum size=0.95cm] {};
\node at (8,0) [circle,draw=purple!100,fill=purple!0,minimum size=0.95cm] {};
\node at (10,0) [circle,draw=purple!100,fill=purple!0,minimum size=0.95cm] {};

\node (a) at (1,1) {1};
\node (b) at (3,1) {4};
\node (c) at (5,1) {16};
\node (d) at (7,1) {68};
\node (e) at (9,1) {304};
\node (f) at (2,2) {1};
\node (g) at (4,2) {6};
\node (h) at (6,2) {30};
\node (i) at (8,2) {146};
\node (j) at (3,3) {1};
\node (k) at (5,3) {8};
\node (l) at (7,3) {48};
\node (m) at (4,4) {1};
\node (n) at (6,4) {10};
\node (o) at (5,5) {1};
\node (p) at (0,0) {1};
\node (q) at (2,0) {2};
\node (r) at (4,0) {6};
\node (s) at (6,0) {22};
\node (t) at (8,0) {90};
\node (u) at (10,0) {394};
\end{tikzpicture}
\begin{figure}[h!]
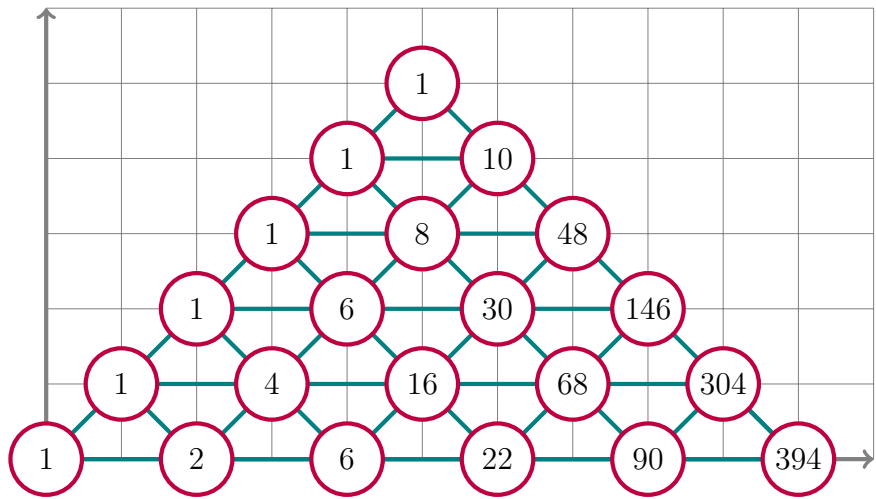

\caption{Large Schr\"{o}der Path, S(x)}
\end{figure}
\end{center}

\begin{center}
\begin{tikzpicture}[ultra thick]
\draw [help lines] (0,0) grid (11,6);
\draw[->] [color=gray] (0,0)--(0,6);
\draw[->] [color=gray] (0,0)--(11,0);
\coordinate (aux) at (0,0);
\draw[->] [color=teal] (0,0)--(5,5);
\draw[->] [color=teal] (5,5)--(10,0);
\draw[->] [color=teal] (1,1)--(2,0);
\draw[->] [color=teal] (2,0)--(3,1);
\draw[->] [color=teal] (3,1)--(4,0);
\draw[->] [color=teal] (4,0)--(5,1);
\draw[->] [color=teal] (5,1)--(6,0);
\draw[->] [color=teal] (6,0)--(7,1);
\draw[->] [color=teal] (7,1)--(8,0);
\draw[->] [color=teal] (8,0)--(9,1);
\draw[->] [color=teal] (9,1)--(10,0);
\draw[->] [color=teal] (2,2)--(3,1);
\draw[->] [color=teal] (3,1)--(4,2);
\draw[->] [color=teal] (4,2)--(5,1);
\draw[->] [color=teal] (5,1)--(6,2);
\draw[->] [color=teal] (6,2)--(7,1);
\draw[->] [color=teal] (7,1)--(8,2);
\draw[->] [color=teal] (3,3)--(4,2);
\draw[->] [color=teal] (4,2)--(5,3);
\draw[->] [color=teal] (5,3)--(6,2);
\draw[->] [color=teal] (6,2)--(7,3);
\draw[->] [color=teal] (4,4)--(5,3);
\draw[->] [color=teal] (5,3)--(6,4);
\draw[->] [color=gray] (2,0)--(4,0);
\draw[->] [color=gray] (4,0)--(6,0);
\draw[->] [color=gray] (6,0)--(8,0);
\draw[->] [color=gray] (8,0)--(10,0);

\draw[->] [color=teal] (1,1)--(9,1);
\draw[->] [color=teal] (2,2)--(8,2);
\draw[->] [color=teal] (3,3)--(7,3);
\draw[->] [color=teal] (4,4)--(6,4);

\node at (1,1) [circle,draw=purple!100,fill=purple!0,minimum size=0.95cm] {};
\node at (3,1) [circle,draw=purple!100,fill=purple!0,minimum size=0.95cm] {};
\node at (5,1) [circle,draw=purple!100,fill=purple!0,minimum size=0.95cm] {};
\node at (7,1) [circle,draw=purple!100,fill=purple!0,minimum size=0.95cm] {};
\node at (9,1) [circle,draw=purple!100,fill=purple!0,minimum size=0.95cm] {};
\node at (2,2) [circle,draw=purple!100,fill=purple!0,minimum size=0.95cm] {};
\node at (4,2) [circle,draw=purple!100,fill=purple!0,minimum size=0.95cm] {};
\node at (6,2) [circle,draw=purple!100,fill=purple!0,minimum size=0.95cm] {};
\node at (8,2) [circle,draw=purple!100,fill=purple!0,minimum size=0.95cm] {};
\node at (3,3) [circle,draw=purple!100,fill=purple!0,minimum size=0.95cm] {};
\node at (5,3) [circle,draw=purple!100,fill=purple!0,minimum size=0.95cm] {};
\node at (7,3) [circle,draw=purple!100,fill=purple!0,minimum size=0.95cm] {};
\node at (4,4) [circle,draw=purple!100,fill=purple!0,minimum size=0.95cm] {};
\node at (6,4) [circle,draw=purple!100,fill=purple!0,minimum size=0.95cm] {};
\node at (5,5) [circle,draw=purple!100,fill=purple!0,minimum size=0.95cm] {};
\node at (0,0) [circle,draw=purple!100,fill=purple!0,minimum size=0.95cm] {};
\node at (2,0) [circle,draw=purple!100,fill=purple!0,minimum size=0.95cm] {};
\node at (4,0) [circle,draw=purple!100,fill=purple!0,minimum size=0.95cm] {};
\node at (6,0) [circle,draw=purple!100,fill=purple!0,minimum size=0.95cm] {};
\node at (8,0) [circle,draw=purple!100,fill=purple!0,minimum size=0.95cm] {};
\node at (10,0) [circle,draw=purple!100,fill=purple!0,minimum size=0.95cm] {};

\node (a) at (1,1) {1};
\node (b) at (3,1) {3};
\node (c) at (5,1) {11};
\node (d) at (7,1) {45};
\node (e) at (9,1) {197};
\node (f) at (2,2) {1};
\node (g) at (4,2) {5};
\node (h) at (6,2) {23};
\node (i) at (8,2) {107};
\node (j) at (3,3) {1};
\node (k) at (5,3) {7};
\node (l) at (7,3) {39};
\node (m) at (4,4) {1};
\node (n) at (6,4) {9};
\node (o) at (5,5) {1};
\node (p) at (0,0) {1};
\node (q) at (2,0) {1};
\node (r) at (4,0) {3};
\node (s) at (6,0) {11};
\node (t) at (8,0) {45};
\node (u) at (10,0) {197};
\end{tikzpicture}
\begin{figure}[h!]
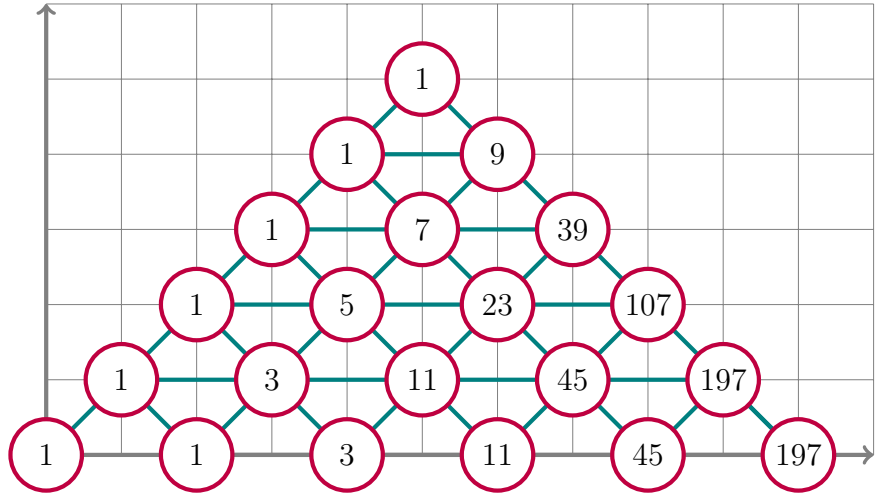

\caption{Small Schr\"{o}der Path, s(x)}
\end{figure}
\end{center}
\noindent The following equations will be proven to be the generating functions for the Schr\"{o}der numbers in Section 2:

\begin{equation}
S(x)=\frac{1-x-\sqrt{1-6x+x^2}}{2x} \label{eq:3}
\end{equation}

\begin{equation}
s(x)=\frac{1+x-\sqrt{1-6x+x^2}}{4x} \label{eq:4}
\end{equation}

\noindent These generate the infinite sequence of the Schr\"{o}der numbers below, which can also be found along the $x$-axis of the lattices above: 

\begin{equation}
S(x)=1+2x+6x^2+22x^3+90x^4+394x^5+...
\end{equation}

\begin{equation}
s(x)=1+x+3x^2+11x^3+45x^4+197x^5+...\\
\end{equation}

\subsection{Motzkin Numbers}

Motzkin paths, denoted $M(x)$ or $M$, are similar in nature to the previous paths, but are composed of the vectors:
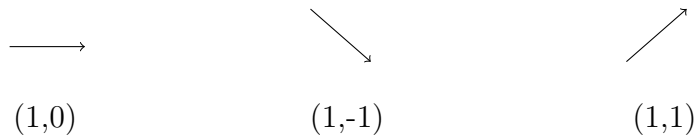
\begin{center}

\begin{tikzpicture}

\draw[->] (0,-.5) -- (1,-.5);

\draw [->](4,0) -- (4.8,-0.7);

\draw [->](8.2,-0.7) -- (9, 0);
\end{tikzpicture}
\end{center}
\tab[3.8cm](1,0)\tab[3.1cm](1,-1)\tab[3.3cm](1,1)\\

\noindent Because the length of the horizontal vector is only half that of the $(2,0)$ level step found in the Schr\"{o}der numbers, there are twice as many of these vectors in Motzkin paths. Together, the above paths combine to form the following lattice in the first quadrant of the coordinate plane from $(0,0)$ to $(n,0)$:\\
\begin{center}
\begin{tikzpicture}[ultra thick]
\draw [help lines] (0,0) grid (9,5);
\draw[->] [color=gray] (0,0)--(0,5);
\draw[->] [color=gray] (0,0)--(9,0);
\draw[->] [color=teal] (0,0)--(4,4);
\draw[->] [color=teal] (4,4)--(8,0);
\draw[->] [color=teal] (1,1)--(2,0);
\draw[->] [color=teal] (2,0)--(3,1);
\draw[->] [color=teal] (3,1)--(4,0);
\draw[->] [color=teal] (4,0)--(5,1);
\draw[->] [color=teal] (5,1)--(6,0);
\draw[->] [color=teal] (6,0)--(7,1);
\draw[->] [color=teal] (7,1)--(8,0);
\draw[->] [color=teal] (2,2)--(3,1);
\draw[->] [color=teal] (3,1)--(4,2);
\draw[->] [color=teal] (4,2)--(5,1);
\draw[->] [color=teal] (5,1)--(6,2);
\draw[->] [color=teal] (6,2)--(7,1);
\draw[->] [color=teal] (3,3)--(4,2);
\draw[->] [color=teal] (4,2)--(5,3);
\draw[->] [color=teal] (5,3)--(6,2);
\draw[->] [color=teal] (4,4)--(5,3);
\draw[->] [color=teal] (0,0)--(9,0);
\draw[->] [color=teal] (1,0)--(2,1);
\draw[->] [color=teal] (2,1)--(3,2);
\draw[->] [color=teal] (3,2)--(4,3);
\draw[->] [color=teal] (4,3)--(5,4);
\draw[->] [color=teal] (2,1)--(3,0);
\draw[->] [color=teal] (3,0)--(4,1);
\draw[->] [color=teal] (4,1)--(5,0);
\draw[->] [color=teal] (5,0)--(6,1);
\draw[->] [color=teal] (3,2)--(4,1);
\draw[->] [color=teal] (5,1)--(6,2);
\draw[->] [color=teal] (4,1)--(5,2);
\draw[->] [color=teal] (5,2)--(6,3);
\draw[->] [color=teal] (6,1)--(7,0);
\draw[->] [color=teal] (7,0)--(8,1);
\draw[->] [color=teal] (6,1)--(7,2);
\draw[->] [color=teal] (4,3)--(5,2);
\draw[->] [color=teal] (5,2)--(6,1);

\draw[->] [color=teal] (0,0)--(1,0);
\draw[->] [color=teal] (1,0)--(2,0);
\draw[->] [color=teal] (2,0)--(3,0);
\draw[->] [color=teal] (3,0)--(4,0);
\draw[->] [color=teal] (4,0)--(5,0);
\draw[->] [color=teal] (5,0)--(6,0);
\draw[->] [color=teal] (6,0)--(7,0);
\draw[->] [color=teal] (7,0)--(8,0);
\draw[->] [color=teal] (1,1)--(2,1);
\draw[->] [color=teal] (2,1)--(3,1);
\draw[->] [color=teal] (3,1)--(4,1);
\draw[->] [color=teal] (4,1)--(5,1);
\draw[->] [color=teal] (5,1)--(6,1);
\draw[->] [color=teal] (6,1)--(7,1);
\draw[->] [color=teal] (2,2)--(3,2);
\draw[->] [color=teal] (3,2)--(4,2);
\draw[->] [color=teal] (4,2)--(5,2);
\draw[->] [color=teal] (3,3)--(4,3);
\draw[->] [color=teal] (4,3)--(5,3);

\node at (1,1) [circle,draw=purple!100,fill=purple!0,minimum size=0.7cm] {};
\node at (3,1) [circle,draw=purple!100,fill=purple!0,minimum size=0.7cm] {};
\node at (5,1) [circle,draw=purple!100,fill=purple!0,minimum size=0.7cm] {};
\node at (7,1) [circle,draw=purple!100,fill=purple!0,minimum size=0.85cm] {};
\node at (2,2) [circle,draw=purple!100,fill=purple!0,minimum size=0.7cm] {};
\node at (4,2) [circle,draw=purple!100,fill=purple!0,minimum size=0.7cm] {};
\node at (6,2) [circle,draw=purple!100,fill=purple!0,minimum size=0.7cm] {};
\node at (3,3) [circle,draw=purple!100,fill=purple!0,minimum size=0.7cm] {};
\node at (5,3) [circle,draw=purple!100,fill=purple!0,minimum size=0.7cm] {};
\node at (4,4) [circle,draw=purple!100,fill=purple!0,minimum size=0.7cm] {};

\node at (0,0) [circle,draw=purple!100,fill=purple!0,minimum size=0.7cm] {};
\node at (2,0) [circle,draw=purple!100,fill=purple!0,minimum size=0.7cm] {};
\node at (4,0) [circle,draw=purple!100,fill=purple!0,minimum size=0.7cm] {};
\node at (6,0) [circle,draw=purple!100,fill=purple!0,minimum size=0.7cm] {};
\node at (8,0) [circle,draw=purple!100,fill=purple!0,minimum size=0.85cm] {};

\node at (1,0) [circle,draw=purple!100,fill=purple!0,minimum size=0.7cm] {};
\node at (3,0) [circle,draw=purple!100,fill=purple!0,minimum size=0.7cm] {};
\node at (5,0) [circle,draw=purple!100,fill=purple!0,minimum size=0.7cm] {};
\node at (7,0) [circle,draw=purple!100,fill=purple!0,minimum size=0.85cm] {};
\node at (2,1) [circle,draw=purple!100,fill=purple!0,minimum size=0.7cm] {};
\node at (4,1) [circle,draw=purple!100,fill=purple!0,minimum size=0.7cm] {};
\node at (6,1) [circle,draw=purple!100,fill=purple!0,minimum size=0.7cm] {};
\node at (3,2) [circle,draw=purple!100,fill=purple!0,minimum size=0.7cm] {};
\node at (5,2) [circle,draw=purple!100,fill=purple!0,minimum size=0.7cm] {};
\node at (4,3) [circle,draw=purple!100,fill=purple!0,minimum size=0.7cm] {};

\node (a) at (1,1) {1};
\node (b) at (3,1) {5};
\node (c) at (5,1) {30};
\node (d) at (7,1) {196};
\node (e) at (9,1) {};
\node (f) at (2,2) {1};
\node (g) at (4,2) {9};
\node (h) at (6,2) {69};
\node (j) at (3,3) {1};
\node (k) at (5,3) {14};
\node (l) at (7,3) { };
\node (21) at (4,4) {1};

\node (p) at (0,0) {1};
\node (q) at (2,0) {2};
\node (r) at (4,0) {9};
\node (s) at (6,0) {51};
\node (t) at (8,0) {323};
\node (u) at (10,0) {};

\node (1) at (1,0) {1};
\node (2) at (3,0) {4};
\node (3) at (5,0) {21};
\node (4) at (7,0) {127};
\node (5) at (9,0) {};
\node (6) at (2,1) {2};
\node (8) at (4,1) {12};
\node (9) at (6,1) {76};
\node (10) at (8,1) {};
\node (11) at (3,2) {3};
\node (12) at (5,2) {25};
\node (13) at (7,2) {};
\node (14) at (4,3) {4};
\node (15) at (6,3) {};

\end{tikzpicture}
\begin{figure}[h!]
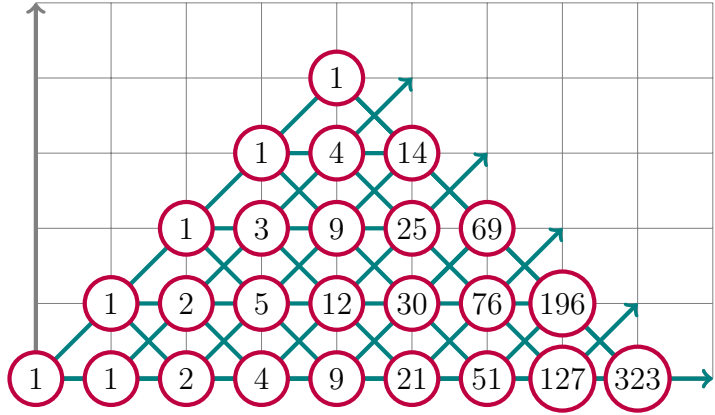

\caption{Motzkin Path, M(x)}
\end{figure}
\end{center}

\noindent The generating function will be proven in Section $2$ to be:

\begin{equation}
M(x)=\frac{1-x-\sqrt{1-2x-3x^2}}{2x^2} \label{eq:7}
\end{equation}

\noindent This generates the following infinite Motzkin sequence, which can also be seen in the numbers along the $x$-axis of the lattice:

\begin{equation}
M(x)=1+x+2x^2+4x^3+9x^4+21x^5+51x^6+...
\end{equation}

\subsection{Pascal's Triangle}

\noindent Pascal's Triangle is an infinite symmetric array that represents the binomial coefficients. To build the triangle, we begin with row zero, which consists of a single one. To build the triangle, each subsequent entry of the next row, $k$, is found by the formula below:

{\normalsize{$$\binom{n-1}{k-1}+\binom{n-1}{k}=\binom{n}{k} = \frac{n!}{k!(n-k)!} $$}}
\\
\noindent Each entry in the triangle can be identified by its position, such as the $n$th row and $k$th column; an entry can be noted as $\binom {n}{k}$. For instance, $\binom{0}{0}=1$ because the value 1 is in row $n=0$ and column $k=0$.

\begin{center}
\def\x{\hspace{5ex}}    
\def\y{\hspace{4.45ex}}  
\def\z{\hspace{3.9ex}}    
\large{\Longstack{
1\\
1\x 1\\
1\x 2\x 1\\
1\x 3\x 3\x 1\\
1\x 4\x 6\x 4\x 1\\
1\x 5\y 10\z 10\y 5\x 1\\
1\x 6\y 15\z 20\z 15\y 6\x 1\\
\vdots
}}
\begin{figure}[h!]
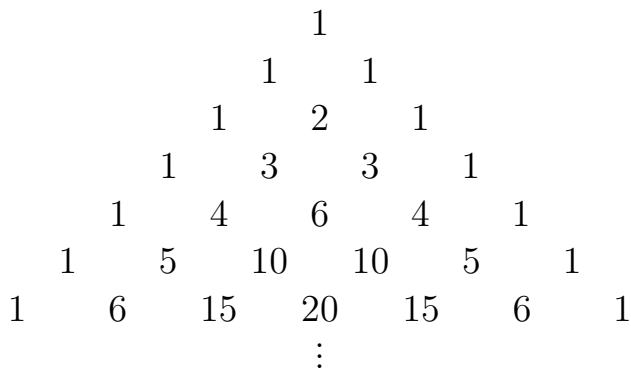

\caption{Pascal's Triangle}
\end{figure}
\end{center}

\section{First Returns \& Generating Functions}
\noindent By definition, a return is the immediate next time the path touches the $x$-axis after leaving it. However, a \textit{complete} return is a single cycle starting at $(0,0)$ that ends on the $x$-axis at $(2n,0)$, with any number of returns happening in between. These returns will be observed in the following subsections. 
\\
\\
\noindent Note: The $n$ represents the number of up steps in a path. Since these paths must remain in the first quadrant of the Cartesian plane, the amount of down steps must equal the amount of up steps. Thus, $2n$ represents the total number of steps.
\subsection{Catalan Numbers}
\noindent The generating function for the Catalan numbers, $C(x)$, can be developed by considering the possible steps the path may take from the origin. When discussing the upward and downward movements that lead us to different values, it must be seen as as a step from $x^n$ to $x^{n+1}$. This can be visualized by following the movements along the Dyck path as seen through the following vectors and in Figure 6 below:\\

\begin{center}
\begin{tikzpicture}[ultra thick]
\draw[->] [color=gray] (0,0)--(0,5);
\draw[->] [color=gray] (0,0)--(10,0);
\draw[->] (0,0)--(1,1) node[midway,above]{\Large{$x$} \ };
\draw [nearly transparent, dashed] (1,1) -- (5,1);
\draw[color=teal] (1,1) sin (1.5,2) cos (2,1.5) sin (2.5,3) cos (3,2) sin (3.5,3.5) cos (4,2) sin (4.5,3) cos (5,1);
\draw[->] (5,1)--(6,0);
\draw[color=teal] (6,0) sin (6.5,1) cos (7,0.5) sin (7.5,1.75) cos (8,1) sin (8.5,1.5) cos (9,1) sin (9.5,0.5) cos (10,0);
\draw[-] (3,3.4)--(3,3.4) node[midway, above]{\Large{$C(x)$}};
\draw[-] (8,1.8)--(8,1.8) node[midway, above]{\Large{$C(x)$}};

\end{tikzpicture}

\begin{figure}[h!]
\caption{First Returns of C(x)}
\end{figure}
\end{center}

\noindent It can be seen that the movements from the origin will either be none, in which case the resulting value is $1$, or the movement will be along the $(1,1)$ vector. A movement along the $(1,1)$ vector causes us to attain an $x$, and once that movement is made, then there are all of the possibilities of $C(x)$. After the path returns to the $x$-axis, it continues on any random path of $C(x)$ and makes any number of returns. This makes the second possibility have a value of $xC(x)\cdot C(x)$. Combining these two options, the first full return is notated as the following: 
\\
$$C(x)=1+xC(x)\cdot C(x)$$

\begin{equation}
C(x)=1+xC^2(x)
\end{equation}

\noindent 
\begin{proof} Consider the figure above. Using the pattern produced by returns, we will define a formula for $C(x)$.
\begin{align*}
C(x)&=1+xC(x)+x^2C^2(x)+x^3C^3(x)+...\\
\\
 &=1+xC(x)+\Big(xC(x)\Big)^2+\Big(xC(x)\Big)^3+...
\end{align*}

\begin{center}
$C(x)=$ $\displaystyle \sum_{n=0}^{\infty} \Big(xC(x)\Big)^n$\\
\vspace{5mm}
$C(x) =$ {\scalebox{1.5}{$\frac{1}{1-xC(x)}$}}\\
\vspace{5mm}
$C(x) = 1+xC^2(x)$\\
\vspace{5mm}
$xC^2(x)-C(x)+1 = 0$\\
\end{center}
\noindent The generating function follows from the quadratic formula and yields Equation 1 as seen in Section 1:
\begin{equation}
C(x)=\frac{1-\sqrt{1-4x}}{2x}\tag{\ref{eq:1}}
\end{equation}
\end{proof}

\subsection{Schr\"{o}der Numbers}
\noindent The generating function for the Schr\"{o}der numbers, $S(x)$ and $s(x)$, can also be developed by considering the possible steps the path may take from the origin. The equation will differ slightly from the Catalan numbers, as we have to consider a different return pattern for both large and small Schr\"{o}der numbers. The upward paths for both large and small Schr\"{o}der numbers should similarly be taken as a step from $x^n$ to $x^{n+1}$, as should a horizontal movement along the $(2,0)$ vector. This can also be visualized by following the movements along the Schr\"{o}der path.


\subsubsection{Large Schr\"{o}der Numbers}

\begin{center}
\begin{tikzpicture}[ultra thick]
\draw[->] [color=gray] (0,0)--(0,5);
\draw[->] [color=gray] (0,0)--(10,0);
\draw[->] (0,0)--(1,1) node[midway,above]{\Large{$x$} \ \ \ };
\draw [nearly transparent, dashed] (1,1) -- (5,1);
\draw[->] (0,0)--(2,0);
\draw[color=teal] (1,1) sin (1.5,2) cos (2,1.5) sin (2.5,3) cos (3,2) sin (3.5,3.5) cos (4,2) sin (4.5,3) cos (5,1);
\draw[->] (5,1)--(6,0);
\draw[color=teal] (6,0) sin (6.5,1) cos (7,0.5) sin (7.5,1.75) cos (8,1) sin (8.5,1.5) cos (9,1) sin (9.5,0.5) cos (10,0);
\draw[double, dotted, color=purple] (2,0) sin (2.5,1) cos (3,0.5) sin (3.5,0.75) cos (4,0.1) cos (4.5,0.3) sin (5,1) sin (5.5,.15) cos (6,0.85) sin (6.5,1.4) cos (7,.35) sin (7.5,1) cos (8,0.5) sin (8.5,0.1) cos (9,0.6) sin (9.5,1.5) cos (10,1.2);
\draw[-] (3,3.4)--(3,3.4) node[midway, above]{\Large{$S(x)$}};
\draw[-] (8,1.8)--(8,1.8) node[midway, above]{\Large{$S(x)$}};
\draw[-] (4,1)--(4,1) node[midway, above]{\Large{$S(x)$}};
\draw[->] (0,0)--(2,0) node[midway,below]{\Large{$x$}};

\end{tikzpicture}

\begin{figure}[h!]
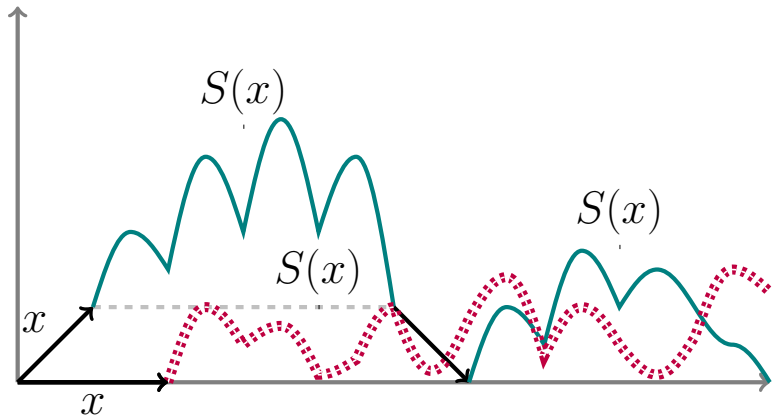

\caption{First Returns of S(x)}
\end{figure}
\end{center}

\begin{proof}
\noindent When analyzing the large Schr\"{o}der numbers, it is evident that the movements from the origin will either be none (in which case the resulting value is 1), horizontally along one or more $(2,0)$ vectors, or upward along the $(1,1)$ vector. A horizontal move would have to attain one $x$ to have the same value as an up and down movement. Once that movement is made, then there are of all the possibilities of $S(x)$. This gives a horizontal movement the total value of $xS(x)$ If the path takes an upward step, it continues on any random path of $S(x)$, then returns to the $x$-axis, and continues again on any random path of $S(x)$ making as many returns as it wants. This gives an upward movement the total value of $xS(x)\cdot S(x)$. Combining these possibilities, the movement of $S(x)$ can be separated by this first movement, thus, the first full return is notated as the following:
$$S(x)=1+xS(x)+xS(x)\cdot S(x)$$
\begin{equation}
S(x)=1+xS(x)+xS^2(x)
\end{equation}
\\


\noindent Use of the quadratic formula yields the generating function as seen in Section 1: 
\begin{equation}
S(x)=\frac{1-x-\sqrt{1-6x+x^2}}{2x}\tag{\ref{eq:3}}
\end{equation}
\end{proof}

\subsubsection{Small Schr\"{o}der Numbers}
\noindent The small Schr\"{o}der numbers have a slightly different pattern than the large Schr\"{o}der numbers because they do not have a level step on the $x$-axis, as it affects the first return and therefore its generating function. 

\begin{center}
\begin{tikzpicture}[ultra thick]
\draw[->] [color=gray] (0,0)--(0,5);
\draw[->] [color=gray] (0,0)--(10,0);
\draw[->] (0,0)--(1,1) node[midway,above]{\Large{$x$} \ \ };
\draw [nearly transparent, dashed] (1,1) -- (5,1);
\draw[color=teal] (1,1) sin (1.5,2) cos (2,1.5) sin (2.5,3) cos (3,2) sin (3.5,3.5) cos (4,2) sin (4.5,3) cos (5,1);
\draw[->] (5,1)--(6,0);
\draw[color=orange] (6,0) sin (6.5,1) cos (7,0.5) sin (7.5,1.75) cos (8,1) sin (8.5,1.5) cos (9,1) sin (9.5,0.5) cos (10,0);
\draw[-] (3,3.4)--(3,3.4) node[midway, above]{\textbf{\Large{$S(x)$}}};
\draw[-] (8,1.8)--(8,1.8) node[midway, above]{{\textbf{\Large{$s(x)$}}}};

\end{tikzpicture}

\begin{figure}[h!]
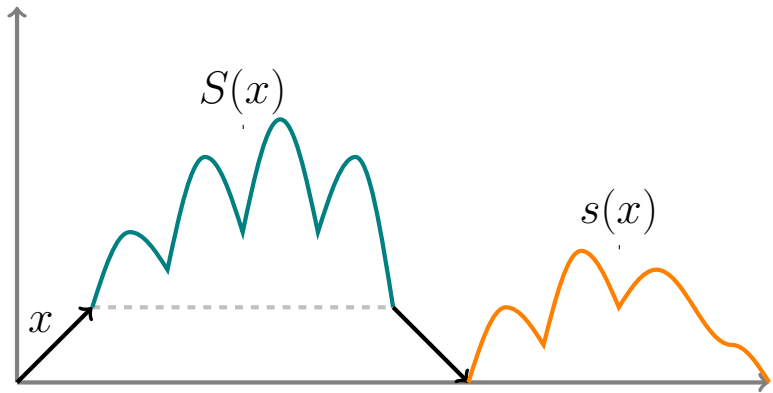

\caption{First Returns of s(x)}
\end{figure}
\end{center}

\noindent When considering the small Schr\"{o}der numbers, it is evident that the movements from the origin will either be none, in which case the resulting value is 1, or the movement will be along the $(1,1)$ vector. There is no initial horizontal move as there are no level steps on the $x$-axis. Once the path takes an upward step, it continues on any random path of $S(x)$ rather than $s(x)$, as there are $(2,0)$ level steps after the path leaves the x-axis. It then returns to the $x$-axis, and continues again on any random path of $s(x)$ making any number of returns. Thus the upward step has the value of $xS(x)\cdot s(x)$. Combining these possibilities, the first full return of $s(x)$ is notated as the following:
\begin{equation}
s(x)=1+xS(x)\cdot s(x)
\end{equation}

\noindent Provided the generating function for large Schr\"{o}der numbers, denoted by $S(x)$, we can then find the generating function for small Schr\"{o}der numbers. \\
\begin{proof}Consider the figure above. Using the pattern produced by returns, we will define a formula for $s(x)$.
\begin{center}
    
$s(x)$ = {\scalebox{1.5}{$\frac{1}{1-xS(x)}$}}
\end{center}
\vspace{.3cm}
\noindent Substituting the generating function of $S(x)$ into the equation above and simplifying, we get the generating function for small Schr\"{o}der numbers as seen in Section 1:

\begin{align*}
s(x)& = {\scalebox{1.65}{$\frac{1}{1-x\Big(\frac{1-x-\sqrt{x^2-6x+1}}{2x}\Big)}$}}\\
\vspace{2.5mm} \\
&= {\scalebox{1.5}{$\frac{2}{1+x+\sqrt{x^2-6x+1}}$}}\\
\vspace{2.5mm}\\
& = {\scalebox{1.5}{$\frac{1+x-\sqrt{x^2-6x+1}}{4x}$}}\tag{\ref{eq:4}}
\end{align*}
\end{proof}

\subsection{Motzkin Numbers}
\noindent The generating function $M(x)$ can also be proven by considering the possible steps it may take from the origin. When we consider the upward and lateral moves that lead to different values, it must be seen as a step from $ x^n$ to $x^{n+1}$. The equation will differ from the Catalan and Schr\"{o}der numbers as the level steps are $(1,0)$. This can also be visualized with the following vectors or using Figure 9.


\begin{center}
\begin{tikzpicture}[ultra thick]
\draw[->] [color=gray] (0,0)--(0,5);
\draw[->] [color=gray] (0,0)--(10,0);
\draw[->] (0,0)--(1,0) node[midway,below]{\Large{$x$}};
\draw[->] (5,1)--(6,0) node[midway,above]{\ \ \Large{$x$}};
\draw [nearly transparent, dashed] (1,1) -- (5,1);
\draw[->] (0,0)--(1,0);
\draw[color=teal] (1,1) sin (1.5,2) cos (2,1.5) sin (2.5,3) cos (3,2) sin (3.5,3.5) cos (4,2) sin (4.5,3) cos (5,1);
\draw[->] (5,1)--(6,0);
\draw[color=teal] (6,0) sin (6.5,1) cos (7,0.5) sin (7.5,1.75) cos (8,1) sin (8.5,1.5) cos (9,1) sin (9.5,0.5) cos (10,0);
\draw[double, dotted, color=green] (1,0) sin (1.5,1) cos (2,0.5) sin (2.5,0.75) cos (3,0.1) cos (3.5,0.3) sin (4,1) cos (4.5,0) sin (5,1) sin (5.5,.15) cos (6,0.85) sin (6.5,1.4) cos (7,.35) sin (7.5,1) cos (8,0.5) sin (8.5,0.1) cos (9,0.6) sin (9.5,1.5) cos (10,1.2);
\draw[-] (3,3.4)--(3,3.4) node[midway, above]{\Large{$M(x)$}};
\draw[-] (8,1.8)--(8,1.8) node[midway, above]{\Large{$M(x)$}};
\draw[-] (3.5,1)--(3.5,1) node[midway, above]{\Large{$M(x)$}};
\draw [->] (0,0)--(1,1) node[midway,above]{\Large{$x$} \ \ };

\end{tikzpicture}

\begin{figure}[h!]
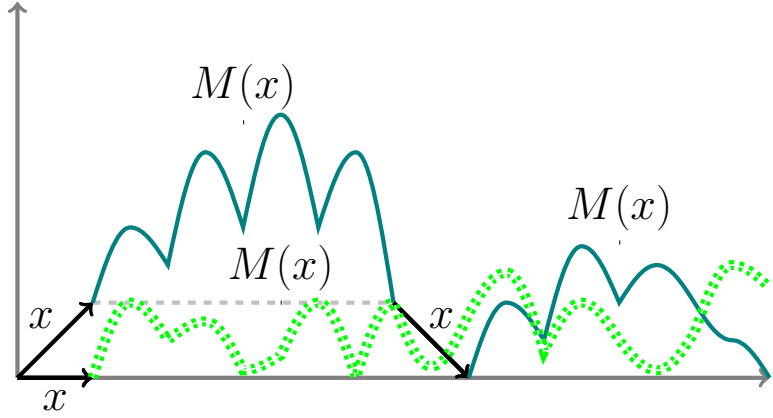

\caption{First Returns of M(x)}
\end{figure}
\end{center}

\begin{proof}
\noindent Thus the movements from the origin will either be none, in which case the resulting value is 1, horizontally along one or more $(1,0)$ vectors, or along the $(1,1)$ vector. A movement along the $(1,1)$ vector gains an $x$ from both the up and down step, and an $M(x)$ where it can do anything. This results in an $x^2M(x)$. Once one movement is made, there are than all the possibilities of $M(x)$, resulting in $x^2M^2(x)$ A vertical initial movement can be represented by $xM(x)$, as discussed in previous sections. A horizontal move would have to attain one $x$ or two $x's$ for the two horizontal movements it needs to make to have the same value as an up and down movement. This yields a total value of $xM(x)$ after an initial horizontal movement. Combining these three possibilities, the movements of $M(x)$ based on the first movement yield the function: 
\begin{equation}
M(x)=1+xM(x)+x^2M^2(x)
\end{equation}

\noindent $M(x)$ can be subtracted from both sides, creating a quadratic function that can be solved using the quadratic formula, yielding Equation 7 as seen in Section 1:

\begin{equation}
M(x)=\frac{1-x-\sqrt{1-2x-3x^2}}{2x^2}\tag{\ref{eq:7}}
\end{equation}
\end{proof}



\section{Recursions}
After analyzing the paths produced in Section 1, we recognized that the terms on the first row along the $x$-axis were the coefficients of Catalan, Schr\"{o}der, and Motzkin sequences, respectively. Upon further investigation, a pattern became apparent where, for example, the second row of coefficients is produced when the sequence from the first row was raised to the second power, and so on. This pattern also held for Pascal's Triangle.\\
\\
A generalized recursion was constructed to determine the coefficients that will appear when any two or more rows (sequences) are multiplied together. Here $G(x)$, also denoted $G$, represents a generalized case. The short notations for Catalan, Schr\"{o}der, and Motzkin numbers will be used in this section, if necessary, to prevent confusion. The recursion works as follows: $$G_k=x^k\cdot G^{k+1}$$
Here, the index begins at 0, meaning the very bottom row is denoted as row zero. Hence, $G_0=G$. Overall, it was found that the subsequent rows are shifted by a power of x; for example, $G_2$ would be $x^2$ multiplied by the sequence at row zero, cubed, i.e. $x^2\cdot (G_0)^3$, which equals $x^2\cdot G^3$. This generalized recursion case will be proven for each set of numbers below.


\subsection{Catalan Numbers}
The primary observation in relation to the numbers generated by the Dyck path, or the number of paths it takes to get to a specific point, can be represented by the recursive equation: 
\begin{equation}
C_k=xC_{k-1}+C_{k+1}
\end{equation}
This states that each value on any path is found by adding the values at the other end of the leftward $(1,1)$ and $(1,-1)$ vectors, not neglecting the $x$ that is attained along the upward vector.
It can then be observed in the first two rows of $C$ that they differ only by 1. We denote the $n^{th}$ row of the Cartesian plane on which the Dyck paths lie by $C_n$. Thus, this observation yields the equation:
\begin{equation}
C_0=1+C_1
\end{equation}
Equation 9 can be solved for $xC^2$ to yield
$$xC^2 = C - 1$$
which will be used alongside Equation 13 and 14 to prove that 
\begin{equation}
C_k = x^k\cdot \big(C_0\big)^{k+1}
\end{equation}

\begin{proof}
\noindent Let $C_k$ represent the horizontal rows of the first quadrant of the coordinate plane. We want to show that equation 15 is true for all positive integers, $k$. \\
\\Base Case: When $k=0$, it is true that $C_0=x^0\cdot \big(C_0\big)^{1}$. \\
\\Thus, we can move on to the inductive step of the proof to show that equation 15 is true for $k+1$. \\
\\\underline{Inductive Step} \\
Assume: 
\begin{align*}
C_{k-1}=x^{k-1}\cdot \big(C_0\big)^k\\
C_{k}=x^{k}\cdot \big(C_0\big)^{k+1}
\end{align*}
Need:
\begin{align*}
C_{k+1}=x^{k+1}\cdot \big(C_0\big)^{k+2}
\end{align*}

\noindent Solving equation 13 for $C_{k+1}$ yields: \\
$$C_{k+1}=C_k-xC_{k-1}$$
From Equation 15 and our assumption, we can replace the right-hand side with:
\begin{align*}
C_{k+1}=x^k\cdot \big(C_0\big)^{k+1}-x\cdot \Big[x^{k-1}\big(C_0\big)^k\Big]\\
C_{k+1}=x^k\cdot \big(C_0\big)^{k+1}-x^k\cdot \big(C_0\big)^k
\end{align*}
Factoring, we get
\begin{align*}
C_{k+1}=x^k\big(C_0\big)^k\cdot (C_0-1)
\end{align*}
From Equation 9, we know that: \\
\begin{align*}
x\big(C_0\big)^2=C_0-1
\end{align*}
Therefore,
\begin{align*}
C_{k+1}=x^k\big(C_0\big)^k\cdot \Big[x\big(C_0\big)^2\Big]\\
C_{k+1}=x^{k+1}\cdot \big(C_0\big)^{k+2}\\
C_{k+1}=x^{k+1}\cdot C^{k+2}
\end{align*}
 Thus, Equation 15 holds true for all integers $k\geq0$.
\end{proof}

\subsection{Schr\"{o}der Numbers}
The primary observation in relation to the numbers generated by the large Schr\"{o}der path, or the number of paths it takes to get to a point, can be represented with the recursive equation: 
\begin{equation}
S_k(x)=S_k=xS_{k-1}+xS_k+S_{k+1}
\end{equation}
This means that each value on any path is found by adding the values at the other end of the leftward $(1,1)$, $(1,-1)$, and $(2,0)$ vectors, not neglecting the x that is attained along the upward and horizontal vector. As the small Schr\"{o}der numbers do not have the level steps on the $x$-axis, this recursive pattern does not work for them in the same way. So we will continue with the large Schr\"{o}der numbers only for this proof. We will begin with Equations 10 and 16 to prove that a similar pattern of recursion exists for Schr\"{o}der numbers: 
\begin{equation}
S_k = x^k\cdot \big(S_0\big)^{k+1}
\end{equation}

\begin{proof}
\noindent Let $S_k$ represent the horizontal rows of the first quadrant of the coordinate plane. We want to show that equation 17 is true for all positive integers, $k$. \\
\\Base Case: When $k=0$, it is true that $S_0=x^0\cdot \big(S_0\big)^{1}$. \\
\\Thus, we can move on to the inductive step of the proof to show that equation 17 is true for $k+1$. \\
\\\underline{Inductive Step} \\
Assume: 
\begin{align*}
S_{k-1}=x^{k-1}\cdot \big(S_0\big)^k\\
S_{k}=x^{k}\cdot \big(S_0\big)^{k+1}
\end{align*}
Need:
\begin{align*}
S_{k+1}=x^{k+1}\cdot \big(S_0\big)^{k+2}
\end{align*}

\noindent Solving equation 16 for $S_{k+1}$ yields: \\
$$S_{k+1}=S_k-xS_k-xS_{k-1}$$
From Equation 17 and the assumption, the right-hand side can be replaced with:
\begin{align*}
S_{k+1}=x^k\cdot \big(S_0\big)^{k+1}-x\cdot \Big[x^k\big(S_0\big)^{k+1}\Big]-x\cdot \Big[x^{k-1}\big(S_0\big)^k\Big]\\
S_{k+1}=x^k\cdot \big(S_0\big)^{k+1}-x^{k+1}\cdot \big(S_0\big)^{k+1}-x^k\cdot \big(S_0\big)^k
\end{align*}
Factoring yields
\begin{align*}
S_{k+1}=x^k\big(S_0\big)^k\cdot (S_0-xS_0-1)
\end{align*}
From Equation 10, it is known that: \\
\begin{align*}
x\big(S_0\big)^2=S_0-xS_0-1
\end{align*}
Therefore,
\begin{align*}
S_{k+1}=x^k\big(S_0\big)^k\cdot \Big[x\big(S_0\big)^2\Big]\\
S_{k+1}=x^{k+1}\cdot \big(S_0\big)^{k+2}\\
S_{k+1}=x^{k+1}\cdot S^{k+2}
\end{align*}
 Thus, Equation 17 holds true for all integers $k\geq0$.
\end{proof}

\subsection{Motzkin Numbers}
The primary observation of the numbers generated by the Motzkin path, or the number of paths it takes to get to a point, can be represented with the recursive equation: 
\begin{equation}
M_k=xM_{k-1}+xM_{k}+xM_{k+1}
\end{equation}
This represents each value on any path that is found by adding the values at the other end of the leftward $(1,1)$, $(1,-1)$, and $(1,0)$ vectors, not neglecting the x that is attained along the upward, downward, and horizontal vectors. We will use Equation 12 to prove that a similar pattern of recursion exists for Motzkin numbers: 
\begin{equation}
M_k = x^k\cdot \big(M_0\big)^{k+1}
\end{equation}

\begin{proof}
\noindent Let $M_k$ represent the horizontal rows of the first quadrant of the coordinate plane. We want to show that equation  is true for all positive integers, $k$. \\
\\Base Case: When $k=0$, it is true that $M_0=x^0\cdot \big(M_0\big)^{1}$ \\
\\Thus, we can move on to the inductive step of the proof to show that equation 19 is true for $k+1$. \\
\\\underline{Inductive Step} \\
Assume: 
\begin{align*}
M_{k-1}=x^{k-1}\cdot \big(M_0\big)^k\\
M_{k}=x^{k}\cdot \big(M_0\big)^{k+1}
\end{align*}
Need:
\begin{align*}
M_{k+1}=x^{k+1}\cdot \big(M_0\big)^{k+2}
\end{align*}

\noindent Solving equation 18 for $xM_{k+1}$ yields: \\
$$xM_{k+1}=M_k-xM_{k}-xM_{k-1}$$
From Equation 19 and the assumption, the right-hand side of the equation can be replaced with:
\begin{align*}
xM_{k+1}=x^k\cdot \big(M_0\big)^{k+1}-x\cdot \Big[x^k\big(M_0\big)^{k+1}\Big]-x\cdot \Big[x^{k-1}\big(M_0\big)^k\Big]\\
xM_{k+1}=x^k\cdot \big(M_0\big)^{k+1}-x^{k+1}\cdot \big(M_0\big)^{k+1}-x^k\cdot \big(M_0\big)^k
\end{align*}
Factoring yields
\begin{align*}
xM_{k+1}=x^k\big(M_0\big)^k\cdot \Big[M_0-xM_0-1\Big]
\end{align*}
From Equation 12, it is known that: \\
\begin{align*}
x^2\big(M_0\big)^2=M_0-xM_0-1
\end{align*}
Therefore,
\begin{align*}
xM_{k+1}=x^k\big(M_0\big)^k\cdot \Big[x^2\big(M_0\big)^2\Big]\\
xM_{k+1}=x^{k+2}\cdot \big(M_0\big)^{k+2}
\end{align*}
Dividing by $x$ on both sides gives:
\begin{align*}
M^{k+1}=x^{k+1}\cdot \big(M_0\big)^{k+2}\\
M^{k+1}=x^{k+1}\cdot M^{k+2}
\end{align*}
Thus, Equation 19 holds true for all integers $k\geq0$.
\end{proof}

\noindent Now, the generating functions and recursive patterns for Catalan, Motzkin and Schr\"{o}der numbers can be combined to describe the relationship between these sequences and Pascal's Triangle. 

\section{Riordan Matrices}
A Riordan matrix is an infinite lower triangular matrix, where the first column is comprised of a function $g(x)$ and each subsequent column has generating function $g(x)\cdot f^k(x)$, where $k$ is the column number, starting with index 0. Pascal's Triangle written as a lower triangular matrix is a Riordan matrix, with the respective generating function:
$$g(x)={\frac{1}{1-x}}$$ for the first column and $$\big(g(x), f(x)\big)=\Big(\frac{1}{1-x}\,,\,\frac{x}{1-x}\Big)$$ for each of the following columns. \\
\begin{center}
$P$ = $\left[ 
\begin{array}{ccccccc}
1 & 0 & 0 & 0 & 0 & 0 &  \\ 
1 & 1 & 0 & 0 & 0 & 0 &  \\ 
1 & 2 & 1 & 0 & 0 & 0 &  \\ 
1 & 3 & 3 & 1 & 0 & 0 & \cdots  \\ 
1 & 4 & 6 & 4 & 1 & 0 &  \\ 
1 & 5 & 10 & 10 & 5 & 1 &  \\ 
&  &  & \cdots  &  &  & 
\end{array}%
\right] $
\end{center}

\noindent A relationship between Pascal's Triangle, denoted $P$, and the large Schr\"{o}der and Motzkin numbers can be seen with respect to their level steps. Since Catalan numbers do not have level steps and small Schr\"{o}der numbers do not follow the recursive pattern, they will not be discussed in this section. Large Schr\"{o}der and Motzkin numbers can be adjusted to have $n$-level steps where for each horizontal step, there are $n$ choices on how to get to the next point. This will be denoted at $S_n$ and $M_n$, respectively. This changes the outcome of the number of paths possible to get to a certain point as instead of there being one level step, there can be $n$ amount of level steps. A new relationship arises out of this, where:

\begin{equation}
PM_n=M_{n+1} \label{Eq20}
\end{equation}

\noindent and, similarly, if Pascal's Triangle is raised to the $n^{th}$ power in matrix form:

\begin{equation}
P^nM=M_{n+1} \label{Eq21}
\end{equation}

\noindent This pattern also holds for $S_n$. These four equations will be proven throughout this section.

\subsection{Step 1: n-colored Level Steps}
In order to define the relationship found between Pascal's Triangle and the n-colored level steps for both Schr\"{o}der and Motzkin numbers, the generating functions for $S_n$ and $M_n$ must first be found.

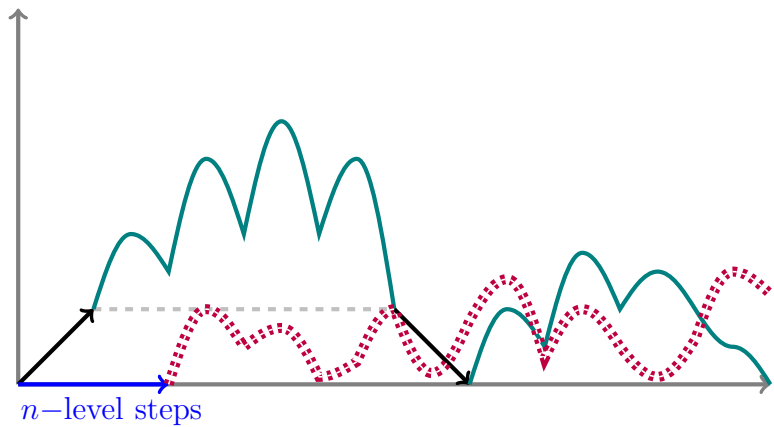
\begin{center}
\begin{tikzpicture}[ultra thick]
\draw[->] [color=gray] (0,0)--(0,5);
\draw[->] [color=gray] (0,0)--(10,0);
\coordinate (aux) at (0,0);
\foreach \i in {1}, 
\draw[->] (aux)--++(1,\i) coordinate (aux);
\draw [nearly transparent, dashed] (1,1) -- (5,1);
\draw[->] [color=blue] (0,0)--(2,0) node[below]{$n-$level steps \ \ \ \ \ \ \ \ \ \ \ };

\draw[color=teal] (1,1) sin (1.5,2) cos (2,1.5) sin (2.5,3) cos (3,2) sin (3.5,3.5) cos (4,2) sin (4.5,3) cos (5,1);
\draw[->] (5,1)--(6,0);
\draw[color=teal] (6,0) sin (6.5,1) cos (7,0.5) sin (7.5,1.75) cos (8,1) sin (8.5,1.5) cos (9,1) sin (9.5,0.5) cos (10,0);
\draw[double, dotted, color=purple] (2,0) sin (2.5,1) cos (3,0.5) sin (3.5,0.75) cos (4,0.1) cos (4.5,0.3) sin (5,1) sin (5.5,.15) cos (6,0.85) sin (6.5,1.4) cos (7,.35) sin (7.5,1) cos (8,0.5) sin (8.5,0.1) cos (9,0.6) sin (9.5,1.5) cos (10,1.2);

\end{tikzpicture}
\begin{figure}
\caption{First Returns of $S_n$}
\end{figure}
\end{center}


\begin{proof}
\noindent Analyzing the large Schr\"{o}der numbers with $n$-level steps, it can be seen that the movements from the origin will either be none, in which case the resulting value is 1, or the movement will be horizontal across the $n$ possible $(2,0)$ vectors, or along the $(1,1)$ vector. A horizontal move would have to attain one $x$ to have the same value as an up and down movement and because there are $n$-level steps, the function also attains an $n$. Once that movement is made, there are then all the possibilities of $S_n(x)$. This yields a total value of $nxS_n(x)$ for an initial horizontal movement. If the path first takes an upward level step, it follows the same pattern it would have made as the original large Schr\"{o}der numbers. Combining these possibilities, the movement of $S_n(x)$ can be separated by this first movement, thus, the first full return is notated as the following: 
$$S_n(x)=1+nxS_n(x)+xS_n(x)\cdot S_n(x)$$
\begin{equation}
S_n(x)=1+nxS_n(x)+xS_n^2(x) \label{eq22}
\end{equation}


\noindent Using the quadratic formula yields the generating function for n-colored level steps for Schr\"{o}der numbers: 
\begin{equation}
S_n(x)=\frac{1-nx-\sqrt{(nx)^2-(2n+4)x+1}}{2x} \label{eq23}
\end{equation}
\end{proof}


\noindent The generating function $M_n(x)$ can also be created by considering the possible steps it may take from the origin. When we consider the upward and lateral moves that lead to different values, it must be seen as a step from $ x^n$ to $x^{n+1}$.

\begin{center}
\begin{tikzpicture}[ultra thick]
\draw[->] [color=gray] (0,0)--(0,5);
\draw[->] [color=gray] (0,0)--(10,0);
\coordinate (aux) at (0,0);
\foreach \i in {1}, 
\draw[->] (aux)--++(1,\i) coordinate (aux);
\draw [nearly transparent, dashed] (1,1) -- (5,1);
\draw[->] [color=orange] (0,0)--(1,0) node[below]{$n-$level steps};
\draw[color=teal] (1,1) sin (1.5,2) cos (2,1.5) sin (2.5,3) cos (3,2) sin (3.5,3.5) cos (4,2) sin (4.5,3) cos (5,1);
\draw[->] (5,1)--(6,0);
\draw[color=teal] (6,0) sin (6.5,1) cos (7,0.5) sin (7.5,1.75) cos (8,1) sin (8.5,1.5) cos (9,1) sin (9.5,0.5) cos (10,0);
\draw[double, dotted, color=green] (1,0) sin (1.5,1) cos (2,0.5) sin (2.5,0.75) cos (3,0.1) cos (3.5,0.3) sin (4,1) cos (4.5,0) sin (5,1) sin (5.5,.15) cos (6,0.85) sin (6.5,1.4) cos (7,.35) sin (7.5,1) cos (8,0.5) sin (8.5,0.1) cos (9,0.6) sin (9.5,1.5) cos (10,1.2);
\end{tikzpicture}

\begin{figure}[h!]
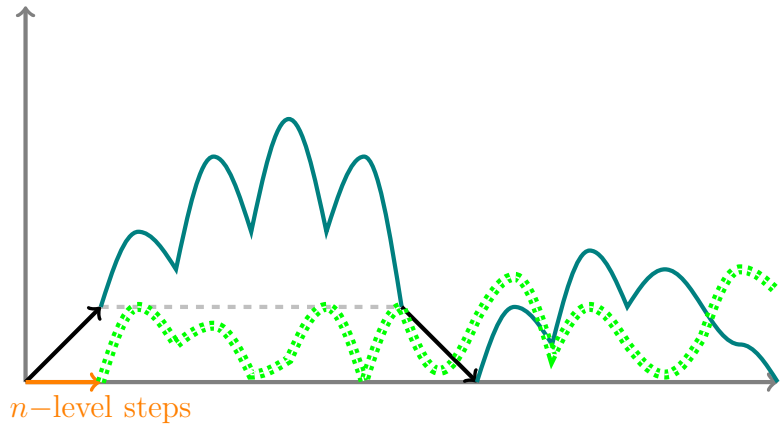

\caption{First Returns of $M_n$}
\end{figure}
\end{center}

\begin{proof}
\noindent This will follow a similar pattern to the original Motzkin numbers, save that the horizontal movement will also gain an $n$ for the $n$-colored level steps. Once one movement is made, there are than all the possibilities of $M_n(x)$. However, a horizontal movement will be able to continue on any random path of $M_n(x)$ at both $(1,0)$ and $(2,0)$.The next possibility is a vertical movement up the $(1,1)$ vector which will pick up an $x^2$ to make up for both the up and down step and then it will continue on a random path before coming back to the $x$-axis and continuing on another random path resulting in $x^2M^2(x)$. Combining these three possibilities, the movements of $M_n(x)$ can be separated by this first movement, giving the function: 

\begin{equation}
M_n(x)=1+nxM(x)+x^2M^2(x) \label{eq24}
\end{equation}

\noindent $M_n(x)$ can be subtracted form both sides, creating a quadratic function that can be solved using the quadratic formula, yielding:

\begin{equation}
M_n(x)=\frac{1-nx-\sqrt{x^2(n^2-4)-2nx+1}}{2x^2} \label{eq25}
\end{equation}
\end{proof}

\noindent Now that $S_n(x)$ and $M_n(x)$ have been proven, they must also be proven for all cases, including $n+1$, which will be shown in the next section. 

\subsection{Step 2: Inductive Reasoning for n$+$1 Level Steps}
The next step in proving the relationship between Pascal's Triangle and Schr\"{o}der and Motzkin numbers is to prove the $n$-level steps case for $n+1$ level steps. The purpose behind this is to show that Equations \ref{eq23} and \ref{eq25} hold for any number of level steps so that there is proof for what the $n+1$ generating function should be when $S_n$ and $M_n$ are multiplied by Pascal's Triangle. Below is the Inductive Proof for Schr\"{o}der Numbers:\\

\begin{proof}
For the following proof, let $S_n=S_n(x)$.
\\
Beginning with Equation 22:
$$S_n(x) = xS_n(x)^2+nxS_n(x)+1$$
Recall that the use of the quadratic formula yields Equation 23:
$$S_n(x)=\frac{1-nx-\sqrt{(nx)^2-x(2n+4)+1}}{2x}$$
We know the following equation to be true due to the recursive pattern that has been proven to exist for the Schr\"{o}der Numbers:
$$S_n = S_{n+1}+xS_{n-1}+nxS_n$$
So when $S_n=S_0$, the previous formula would imply that $S_0(x)=1$.
Continuing with the case n=0 yields the equation
$$xS_n(x)^2+nxS_n(x)+1=1+nxS_n+S_{n+1}$$
Two terms cancel, leaving only $$xS_n(x)^2=S_{n+1}$$
Thus, $S_1=xS_0(x)^2$ and the base case is proven.
\\
\\
\underline{Inductive Step}
\\
Assume:
$$S_n(x)=x^nS^{n+1}$$
$$S_{n-1}(x)=x^{n-1}S^n$$
Need:
\begin{equation}
    S_{n+1}(x)=x^{n+1}S^{n+2}
\end{equation}
\\
\\
The equation 
$$S_n = S_{n+1}+xS_{n-1}+nxS_n$$
can be solved for $S_{n+1}$ and factored by $S_n$ to yield $$S_{n+1}=S_n(1-nx)-xS_{n-1}$$
Using the initial assumptions, this becomes 
$$S_{n+1}=x^nS^{n+1}\big(1-nx\big)-x\big(x^{n-1}S^n\big)$$
Distribution leads to: 
$$S_{n+1}=x^nS^{n+1}-nx^{n+1}S^{n+1}-\big(x^{n}S^n\big)$$
This can be factored, yielding 
$$S_{n+1}=x^nS^n\big(S-nxS-1\big)$$
Equation 22 can be used to replace the parenthesis terms with $xS_n^2$
so that $$S_{n+1}=x^nS^n(xS^2_n)$$
Rearranging the exponents yields the desired result, proving Equation 26 to be true based on the initial assumptions.
$$S_{n+1}=X^{n+1}S_n^{n+1}$$
\end{proof}

\noindent Below is the inductive proof for the Motzkin numbers. 
\begin{proof}
For the following proof, let $M_n=M_n(x)$
\\
Beginning with Equation 24:
$$M_n(x) = x^2M_n(x)^2+nxM_n(x)+1$$
Recall that the use of the quadratic formula yields Equation 25:
$$M_n(x)=\frac{1-nx-\sqrt{x^2(n^2-4)-2nx+1}}{2x^2}$$
We know the following equation to be true due to the recursive pattern that has been proven to exist for the Motzkin Numbers:
$$M_n=xM_{n+1}+xM_{n-1}+nxM_n$$
So when $M_n=M_0$ the previous formula would imply that $M_0(x)=1+nxM_0+xM_1$
This can be set equal to Equation 24
$$1+nxM_0+xM_1=x^2M_n^2+nxM_n+1$$
Two terms cancel, leaving only $$M_1=xM_n^2$$
Thus, the base case is proven.
\\
\\
\underline{Inductive Step}
\\
Assume:
$$M_n=x^nM^{n+1}$$
$$M_{n-1}=x^{n-1}M^n$$
Need:
\begin{equation}
    M_{n+1}(x)=x^{n+1}M^{n+2}
\end{equation}
\\
\\
The equation 
$$M_n = xM_{n+1}+xM_{n-1}+nxM_n$$
can be solved for $xM_{n+1}$ to yield 
$$xM_{n+1}=xM_{n+1}+xM_{n-1}+mxM_n$$
Using the initial assumptions to replace $M_n$ and $M_{n-1}$, this becomes 
$$xM_{n+1}=xM^{n+1}-nx\big(x^nM^{n+1}\big)-x\big(x^{n-1}M^n\big)$$
Distribution leads to: 
$$xM_{n+1}=x^n\big(M^{n+1}-mxM^{n+1}-M^n\big)$$
This can be factored, yielding 
$$xM_{n+1}=x^nM^n\big(M-nxM-1\big)$$
Equation 24 can be used to replace the parenthesis terms with $x^2M^2$
so that $$xM_{n+1}=x^nM^n\big(x^2M^2\big)$$
Both sides can be divided by $x$, giving the desired result and proving Equation 27 to be true based on the initial assumptions.
$$M_{n+1}=x^{n+1}M^{n+2}$$

\end{proof}

\subsection{Step 3: Proving Pascal's Triangle to the nth Power}
In order to continue proving the relationship discovered above, it is necessary to understand what Pascal's Triangle to the $nth$ power would be, and prove it for all cases, including $n+1$. Below is the inductive proof for Pascal's Triangle with $n+1$ rows. 
\begin{proof}
First, we begin with Pascal's Triangle written as a generating function with multiplicative factors.
$$P=\bigg(\frac{1}{1-x}\ ,\ \frac{x}{1-x}\bigg)$$
Therefore, for $P^n$, the generating function with multiplicative factors becomes:
$$P^n=\bigg(\frac{1}{1-nx}\ ,\ \frac{x}{1-nx}\bigg)$$
Now, the inductive proof can be done:
\\Base Case: When $n=0$, it is true that  
$P^0 = \begin{bmatrix}
           {1} \\
           {1} \\
           \vdots \\
           {1}
         \end{bmatrix} $\\
\\Thus, we can move on to the inductive step of the proof to show that it is true for $n+1$. \\
\\\underline{Inductive Step} \\
Assume: 
\begin{align*}
P^n=\bigg(\frac{1}{1-nx}\ ,\ \frac{x}{1-nx}\bigg)\\
\end{align*}
Need:
\begin{align*}
P^{n+1}=\bigg(\frac{1}{1-(n+1)x}\ ,\ \frac{x}{1-(n+1)x}\bigg)
\end{align*}
The method to solving this through induction will be to use a Riordan Matrix, which is defined below:
$$\big(g,f\big)*\big(h,l\big)=\Big(g\cdot h(f),l(f)\Big)$$
In terms of Pascal's Triangle:
$$P*P^n=P^{n+1}$$
That is, we want to solve:
$$\bigg(\frac{1}{1-x}\ ,\ \frac{x}{1-x}\bigg)*\bigg(\frac{1}{1-nx}\ ,\ \frac{x}{1-nx}\bigg)$$
From this Riordan Array, the first part $g\cdot h(f)$ can be solved:\\
\begin{center}
$$g\cdot h(f)=\frac{1}{1-x}\cdot \bigg(\frac{1}{1-n\big(\frac{x}{1-x}\big)}\bigg)$$\\
\end{center}
Simplifying yields:
\begin{align*}
g\cdot h(f)=\frac{1}{1-x}\cdot \bigg(\frac{1-x}{1-x-nx}\bigg)
\end{align*}
Further simplification gives the equation: \\
\begin{align*}
g\cdot h(f)=\frac{1}{1-x(n+1)}
\end{align*}
Next, $l(f)$can be calculated:
\begin{align*}
l(f)=\frac{\frac{x}{1-x}}{1-n\big(\frac{x}{1-x}\big)}
\end{align*}
Simplifying that fully yields:
\begin{align*}
l(f)=\frac{x}{1-x(n+1)}
\end{align*}
Written as a Riordan Matrix, it reads:
\begin{align*}
    \Big(g\cdot h(f),l(f)\Big)=\bigg(\frac{1}{1-x(n+1)},\frac{x}{1-x(n+1)}\bigg)
\end{align*}
Thus, $P^{n+1}$ is true based on the previous assumptions.  
\end{proof}

\subsection{Step 4: Proving the Relationship}
In this section, the relationship found between Pascal's Triangle and the Schr\"{o}der and Motzkin Numbers will be proved. This relationship will define what happens when Pascal's Triangle, as a Riordan Matrix, is multiplied by the Schr\"{o}der or Motzkin numbers with $n$-level steps or vice versa, if we raise Pascal's Triangle to $n^{th}$ power and multiply it by the original Schr\"{o}der or Motzkin numbers. First, Equation \ref{Eq20} will be proven to be true. 
$$PM_n=M_{n+1}$$

\begin{proof}
To prove this relationship, the Riordan Matrix mentioned in the previous section will be used:
$$\big(g,f\big)*\big(h,l\big)=\Big(g\cdot h(f),l(f)\Big)$$
In the end, we want $M_{n+1}$ in the form of a Riordan Matrix, where the generating function, and its' multiplicative factor are given. This is why $M_n$ and its recursive property were proven in the previous section, in order to know what this process should yield. \\
\\ Therefore, $M_{n+1}(x)=$ 
{\small{
\begin{align*}
    {{\Bigg(\frac{1-(n+1)x-\sqrt{x^2\Big((n+1)^2-4\Big)-2x(n+1)+1}}{2x^2}\ ,\ \frac{1-(n+1)x-\sqrt{x^2\Big((n+1)^2-4\Big)-2x(n+1)+1}}{2x}\ \Bigg)}}
\end{align*}}}
We will be multiplying:
\begin{align*}
    P* M_n=\Big(\frac{1}{1-x}\ ,\ \frac{x}{1-x}\Big)* \bigg(\frac{1-nx-\sqrt{x^2(n^2-4)-2nx+1}}{2x^2}\ ,\ \frac{1-nx-\sqrt{x^2(n+^2-4)-2nx+1}}{2x}\bigg)
\end{align*}
Now, solving for $g\cdot h(f)$:

$$\frac{1}{1-x}\cdot \Bigg(\frac{1-n(\frac{x}{1-x})-\sqrt{(\frac{x}{1-x})^2(n^2-4)-2n\big(\frac{x}{1-x}\big)+1}}{2\ \big(\frac{x}{1-x}\big)^2}\Bigg)$$\\
\noindent Simplifying yields:
$$=\bigg(\frac{1}{1-x}\bigg)\cdot \bigg(\frac{1-x-nx-\sqrt{x^2(n^2-4)-2nx(1-x)+(1-x)^2}}{1-x}\bigg)\cdot \bigg(\frac{(1-x)^2}{2x^2}\bigg)$$\\
$$=\frac{1-(n+1)x-\sqrt{x^2(n^2+2n-3)-2x(n+1)+1}}{2x^2}$$
$$=\frac{1-(n+1)x-\sqrt{x^2\Big((n+1)^2-4\Big)-2x(n+1)+1}}{2x^2}$$\\
\noindent This is the desired result for Part 1. Next, $l(f)$ must be proven:\\
$$\frac{1-n(\frac{x}{1-x})-\sqrt{\big(\frac{x}{1-x}\big)^2(n^2-4)-2n\big(\frac{x}{1-x}\big)+1}}{2\ \big(\frac{x}{1-x}\big)}$$\\
\noindent Simplifying, yields:
\begin{center}
$$=\frac{1-x}{2x}\cdot \bigg(\frac{1-x-nx-\sqrt{x^2(n^2-4)-2nx(1-x)+(1-x)^2}}{1-x}\bigg)$$
$$=\frac{1-(n+1)x-\sqrt{x^2(n^2-4+2n+1)-2x(n+1)+1}}{2x}$$
$$=\frac{1-(n+1)x-\sqrt{x^2\Big((n+1)^2-4\Big)-2x(n+1)+1}}{2x}$$
\end{center}

\noindent This is the desired result. Therefore, $P M_n=M_{n+1}$.
\end{proof}
\vspace{.2in}
\noindent Next, the relationship from Equation \ref{Eq21} will be proven and it will be shown that the same answer is produced with $P^nM$.
\begin{proof}
First, we start with:
\begin{align*}
    P^n* M=\bigg(\frac{1}{1-nx}\ ,\ \frac{x}{1-nx}\bigg)* \bigg(\frac{1-x-\sqrt{1-2x-3x^2}}{2x^2}\ ,\  \frac{1-x-\sqrt{1-2x-3x^2}}{2x}\bigg)
\end{align*}

\noindent Next, we use the Riordan Matrix Method to prove $g\cdot h(f)$:
$$=\frac{1}{1-nx}\cdot \Bigg(\frac{1-(\frac{x}{1-nx})-\sqrt{1-2\big(\frac{x}{1-nx}\big)-3\big(\frac{x}{1-nx}\big)^2}}{2\ \big(\frac{x}{1-nx}\big)^2}\Bigg)$$

\noindent Simplifying yields:
$$=\bigg(\frac{1}{1-nx}\bigg)\cdot \bigg(\frac{1-x-nx-\sqrt{(1-nx)^2-2x(1-nx)-3x^2}}{1-nx}\bigg)\cdot \bigg(\frac{(1-nx)^2}{2x^2}\bigg)$$\\
$$=\frac{1-(n+1)x-\sqrt{x^2(n^2+2n-3)-2x(n+1)+1}}{2x^2}$$
$$=\frac{1-(n+1)x-\sqrt{x^2\Big((n+1)^2-4\Big)-2x(n+1)+1}}{2x^2}$$

\noindent Next, we prove $l(f)$:
$$= \frac{1-(\frac{x}{1-nx})-\sqrt{1-2\big(\frac{x}{1-nx}\big)-3\big(\frac{x}{1-nx}\big)^2}}{2\ \big(\frac{x}{1-nx}\big)}$$
Simplifying yields:
$$=\frac{1-nx}{2x}\cdot \frac{1-x-nx-\sqrt{(1-nx)^2-2x(1-nx)-3x^2}}{1-nx}$$
$$=\frac{1-(n+1)x-\sqrt{x^2(n^2+2n-3)-2x(n+1)+1}}{2x}$$
$$=\frac{1-(n+1)x-\sqrt{x^2\Big((n+1)^2-4\Big)-2x(n+1)+1}}{2x}$$
\end{proof}
\noindent This is the desired result. Therefore, $P^n M=M_{n+1}$.\\

\section{Conclusion}
There is a relationship found between different random walks and Pascal's Triangle as can be seen through the last section. These combinatoric patterns have been found through multiple methods in the past, but this paper shows the overarching method that they are all connected to. From first returns and generating functions to Riordan matrices and Pascal's Triangle, we discovered patterns that led us to develop more questions. These numbers and patterns have applications in counting and coding, which is why we wanted to further explore them. \\
\\
Later, we will explore the effect that using multicolored vectors rather than uni-colored vectors will have on Catalan, Schr\"{o}der, and Motzkin sequences, and we will investigate the average returns for these number sequences.

\section{Acknowledgements}
This research is made possible by the generous support of the National Security Agency (NSA), Mathematical Association of America (MAA), National Science Foundation (NSF) grant DMS-1560332 administered though the American Statistical Association (ASA), Delta Kappa Gamma Educational Foundation, and Morgan State University.


\begin{thebibliography}{}
\bibitem{Riordan}\emph{S. Getu, L. Shapiro, W.-J. Woan, L. C. Woodson}, The Riordan group,
Discrete Applied Mathematics Vol. 34.1 (1991), 229-239
\bibitem{peart}\emph{P. Peart, L. Woodson}, Triple factorization of some Riordan matrices, The Fibonacci Quart. Vol. 31 (1339), 121-128
\bibitem{getu}\emph{S. Getu, L. W. Shapiro, L. C. Woodson, W.-J. Woan}, How to guess a generating function, SIAM Journal on Discrete Mathematics Vol. 5.4 (1992)
\bibitem{aigner}\emph{M. Aigner}, Combinatorial Theory, Springer (2013), ISBN 978-3-540-61787-7
\bibitem{stanley1}\emph{R. Stanley}, Enumerative Combinatorics, Cambridge University Press Vol. 2 (1999), ISBN 978-0-521-56069-1
\bibitem{stanley2}\emph{R. Stanley}, Catalan Numbers, Massachusetts Institute of Technology (2015), ISBN 978-1-107-07509-2
\bibitem{oeis}\emph{Sloane's\,Online\,Encyclopedia\,of\,Integer\,Sequences}, http://oeis.org/
\end{thebibliography}
\end{document}